\documentclass[fleqn,preprint,3p,a4paper]{elsarticle}
\usepackage[T1]{fontenc}
\usepackage{lmodern}
\usepackage{microtype}
\usepackage{amsmath,amssymb,amsthm,mathtools,mathrsfs}
\usepackage{enumitem}
\usepackage{booktabs}
\usepackage{xcolor}
\usepackage[colorlinks=true,linkcolor=blue!45!black,citecolor=blue!45!black,urlcolor=blue!45!black]{hyperref}
\hypersetup{pdftitle={H-sober spaces and Hofmann-Mislove type theorems},
  pdfauthor={Xiaoquan Xu and Xinjie Huang},
  pdfsubject={H-sobriety, Smyth power spaces and Hofmann-Mislove representations}}
\biboptions{numbers,sort&compress}
\journal{Topology and its Applications}
\newtheorem{theorem}{Theorem}[section]
\newtheorem{proposition}[theorem]{Proposition}
\newtheorem{lemma}[theorem]{Lemma}
\newtheorem{corollary}[theorem]{Corollary}
\theoremstyle{definition}
\newtheorem{definition}[theorem]{Definition}
\newtheorem{example}[theorem]{Example}
\theoremstyle{remark}
\newtheorem{remark}[theorem]{Remark}
\numberwithin{equation}{section}
\newcommand{\Om}{\mathcal O}
\newcommand{\F}{\mathcal F}
\newcommand{\G}{\mathcal G}
\newcommand{\K}{Q}
\newcommand{\cl}[1]{\overline{#1}}
\newcommand{\N}{\mathbb N}
\newcommand{\PhiX}{\Phi_X}
\newcommand{\Cw}{\mathscr F_\omega}
\newcommand{\OFcb}{\mathsf{OFilt}_{\mathrm{cb}}}
\newcommand{\cb}{\mathrm{co}}
\newcommand{\Hco}{\ensuremath{\mathrm{H}_{\mathrm{co}}}}

\setlist[enumerate]{label=\textup{(\arabic*)},leftmargin=1.7em,itemsep=3pt,topsep=5pt}
\begin{document}

\begin{frontmatter}

\title{H-sober spaces and Hofmann--Mislove type theorems\tnoteref{t1}}
\tnotetext[t1]{This research was supported by the National Natural Science Foundation of China  (Nos. 12471070, 12071199).}

\author[mnnu]{Xiaoquan Xu\corref{cor1}}
\ead{xiqxu2002@163.com}
\author[mnnu]{Xinjie Huang}
\ead{2418782429@qq.com}
\cortext[cor1]{Corresponding author.}
\affiliation[mnnu]{organization={School of Mathematics and Statistics,
  Minnan Normal University},
  city={Zhangzhou},postcode={363000},country={China}}

\begin{abstract}
A $T_0$-space is called $\omega$-sober if every irreducible closed subset
having a countable dense subset is the closure of a unique point.
This is precisely $H$-sobriety for the countable irreducible subset
system $\mathcal R^\omega$. We give two natural Hofmann--Mislove
representations in the general $H$-sobriety framework:
filters induced by $H$-sets of $X$ characterize $H$-sobriety, whereas
filters induced by $H$-sets of the Smyth power space characterize
super $H$-sobriety. Their relation is expressed by property $Q$.
A minimal closed set argument with a countable selection proves that
$\mathcal R^\omega$ has property $Q$. Consequently, $X$ is $\omega$-sober
if and only if its Smyth power space is $\omega$-sober, equivalently,
every countable irreducible compact-family filter has a compact
saturated representation. We also obtain equivalent recovery,
compact-intersection and representation conditions for an enlarged
class containing the proper Scott-open filters with countable open
bases and the filters induced by countable irreducible subsets.
A further equivalent formulation uses traces on countable subspaces.
The filtered-family specializations recover well-filteredness and
$\omega$-well-filteredness. Explicit examples show that
$\omega$-sobriety strictly implies $\omega$-well-filteredness, which
strictly implies representation for countable open bases alone.
\end{abstract}

\begin{keyword}
Irreducible subset system \sep
$H$-sober space \sep
$\omega$-sober space \sep
Hofmann--Mislove theorem \sep
Scott-open filter \sep
Smyth power space
\MSC[2020] 54D10 \sep 54B20 \sep 54D30 \sep 06B35
\end{keyword}

\end{frontmatter}

\section{Introduction}

The Hofmann--Mislove theorem identifies the Scott-open filters of the
open-set lattice of a sober space with its compact saturated subsets
\cite{HM,KP}. Its countable variants require care: countability of a
filter, existence of a countable filter base, and countability of a set
whose closure is irreducible impose different restrictions.

In this paper, for a $T_0$-space $X$, an \emph{open filter} means a Scott-open filter on the
open-set lattice $\Om(X)$. A filter on $\Om(X)$ is \emph{countably generated} if it is generated by
countably many \emph{open members of the filter}. The notion of
$\omega$-sobriety concerns irreducible closed sets with a countable dense
subset, and does
not require those closed sets themselves to be countable. In the language
of irreducible subset systems, this is precisely
$\mathcal R^\omega$-sobriety \cite{XuH}, where
\[
 \mathcal R^\omega(X)
 =\{A\subseteq X:A\text{ is countable, nonempty and irreducible}\}.
\]
Thus $\omega$-sobriety is a special case of $H$-sobriety, with the
countability restriction placed on the inducing irreducible subset.

We distinguish two natural classes of proper Scott-open filters for an
irreducible subset system $H$. The first consists of
$\F_A=\{U\in\Om(X):U\cap A\ne\varnothing\}$ for $A\in H(X)$.
The second consists of
$\PhiX(\mathcal A)=\bigcup_{G\in\mathcal A}\PhiX(G)$ for
$\mathcal A\in H(P_S(X))$, where $P_S(X)$ is the Smyth power space of $X$.
Their representation properties characterize, respectively,
$H$-sobriety of $X$ and $H$-sobriety of $P_S(X)$
(Theorems~\ref{thm:H-point} and~\ref{thm:H-Smyth}).
These are filter formulations of the existing $H$-sobriety and
super $H$-sobriety framework
\cite[Theorems~3.6 and~5.7]{SSXZ},
\cite[Theorem~5.5]{XuH}.
The local property-$Q$ criterion
\cite[Theorem~6.4]{SSXZ} explains when the two conditions coincide.

For an arbitrary $H$, its restriction $H^\omega$ to countable members
is again an irreducible subset system, so both representation theorems
apply. For $\mathcal R^\omega$ we prove property $Q$ directly, using
a minimal closed set and one chosen point from each member of a
countable irreducible family in $P_S(X)$
(Lemma~\ref{lem:Romega-Q}). This gives
\[
 X\text{ is }\omega\text{-sober}
 ~\Longleftrightarrow~
 P_S(X)\text{ is }\omega\text{-sober}.
\]
Theorem~\ref{thm:omega-Smyth} expresses this equivalence in terms of
both point-induced filters and countable irreducible compact-family
filters. In the latter case, \emph{countable} counts the members of the
family, without imposing a countable open base on the induced filter.

Representation for countably generated Scott-open filters alone is
weaker than $\omega$-sobriety. Indeed, it already holds in every
$\omega$-well-filtered space by the theorem of Goubault-Larrecq
\cite[Theorem~1 and the following remark]{GL}. We denote this
open-generator representation property by \Hco, with ``c'' standing for
countable and ``o'' for open. Goubault-Larrecq allows a broader notion of
generation by descending subsets, which we distinguish in
Remark~\ref{rem:GL-scope}. By Theorem~\ref{thm:H-point}, the filters
$\F_A$ induced by countable irreducible subsets already characterize
$\omega$-sobriety. Theorem~\ref{thm:cb} shows that one may enlarge this
test class by adjoining all proper Scott-open filters with countable
open bases. The resulting class $\Cw(X)$ is defined explicitly in
Section~3. The theorem gives equivalent recovery and representation
conditions for this enlarged class, whereas the countably based
filters alone do not characterize $\omega$-sobriety in general.

Theorem~\ref{thm:cb} has a direct minimal-closed-set proof, and
Theorem~\ref{thm:exact} gives a further equivalent formulation using
traces on countable subspaces. The classical, well-filtered and
$\omega$-well-filtered versions are compared through their inducing
families. The examples establish that both implications
$\omega$-sobriety $\Rightarrow$ $\omega$-well-filteredness
$\Rightarrow$ \Hco\ are strict.

The terminology here also differs from \emph{countably sober} in
\cite{YLY, YS}, where the defining irreducibility condition involves
countable closed covers. No identification with that condition is intended.

\section{Preliminaries and Hofmann--Mislove representations}

All spaces are $T_0$, and ``countable'' includes finite.
We write $\N=\{1,2,\ldots\}$. For a space $X$, let $\Om(X)$ denote
its lattice of open subsets, ordered by inclusion, and write $\cl A$
for the closure of $A\subseteq X$.

A nonempty subset $A$ of $X$ is \emph{irreducible} if, whenever
$A\subseteq C_1\cup C_2$ for closed subsets $C_1,C_2$ of $X$,
one has $A\subseteq C_1$ or $A\subseteq C_2$.
Equivalently, for all $U,V\in\Om(X)$,
\[
 A\cap U\ne\varnothing\ \text{and}\ A\cap V\ne\varnothing
 ~\Longrightarrow~ A\cap U\cap V\ne\varnothing.
\]
A subset is irreducible if and only if its closure is irreducible.
A point $x$ is a \emph{generic point} of a closed set $C$ if
$C=\cl{\{x\}}$. The space $X$ is \emph{sober} if every irreducible
closed subset has a unique generic point.

The specialization order is
\[
 x\leq_X y~\Longleftrightarrow~ x\in\cl{\{y\}}.
\]
For $A\subseteq X$, put
\[
 \uparrow A=\{y\in X:a\leq_X y\text{ for some }a\in A\},
 ~ \uparrow x=\uparrow\{x\}.
\]
The \emph{saturation} of $A$ is the intersection of all open sets
containing $A$; it equals $\uparrow A$.
Thus a subset is \emph{saturated} if and only if it is an upper set
for the specialization order. Open sets are upper sets as well.
In particular, if an open or saturated set meets $\cl{\{x\}}$,
then it contains $x$.

A subset $K\subseteq X$ is \emph{compact} if every open cover of $K$
has a finite subcover; no separation assumption is included.
We write $\K(X)$ for the collection of nonempty compact saturated
subsets of $X$.

\begin{definition}\label{def:omega}
A space $X$ is \emph{$\omega$-sober} if every
irreducible closed subset $A$ for which
\[
 A=\cl D~\text{for some countable }D\subseteq X
\]
is $\cl{\{x\}}$ for a unique $x\in X$.
\end{definition}

Equivalently, for every countable irreducible subset $D$ of $X$ there is
a unique $x$ with $\cl D=\cl{\{x\}}$. This follows because a subset is
irreducible if and only if its closure is irreducible. Uniqueness follows
from the $T_0$ assumption. On a countable space, $\omega$-sobriety is
exactly sobriety.
In particular, Definition~\ref{def:omega} is the $H$-sobriety condition
for the countable irreducible subset system $\mathcal R^\omega$;
the general definitions are given below.

A nonempty subset $D$ of a poset is \emph{directed} if every pair
of its members has an upper bound in $D$.
For a family of open sets ordered by inclusion, this means that
any two members are contained in a third member of the family.

A \emph{proper filter} $\F\subseteq\Om(X)$ is a nonempty upper set
closed under finite intersections and not containing $\varnothing$.
It is \emph{Scott-open} if, for every nonempty directed family
$\{U_i:i\in I\}\subseteq\Om(X)$,
\[
 \bigcup_{i\in I}U_i\in\F~\Longrightarrow~
 U_i\in\F\text{ for some }i\in I.
\]
For a subset $K\subseteq X$, put
\[
 \PhiX(K)=\{U\in\Om(X):K\subseteq U\}.
\]

For a nonempty irreducible subset $A\subseteq X$, define
\begin{equation}\label{eq:FA}
 \F_A=\{U\in\Om(X):A\cap U\ne\varnothing\}.
\end{equation}
Irreducibility makes $\F_A$ a proper filter. It is \emph{completely prime}:
if an arbitrary union belongs to $\F_A$, some member of that union
belongs to $\F_A$. In particular, $\F_A$ is Scott-open. If
$\cl A=\cl{\{x\}}$, then
\begin{equation}\label{eq:point-filter}
 \F_A=\{U\in\Om(X):x\in U\}=\PhiX(\uparrow x).
\end{equation}
The set $\uparrow x$ is always nonempty, compact and saturated: any
open set containing $x$ contains all of $\uparrow x$.

\begin{lemma}\label{lem:compact}
If $K\in\K(X)$, then $\PhiX(K)$ is a proper Scott-open filter, and
$K=\bigcap\PhiX(K)$. Conversely, if a proper Scott-open filter $\F$
satisfies
\[
 \F=\PhiX\!\left(\bigcap\F\right),
\]
then $\bigcap\F$ is a nonempty compact saturated set. Thus a compact
saturated set representing $\F$ is unique.
\end{lemma}

\begin{proof}
For a directed open cover of a compact set, a finite subcover is contained
in one member of the directed family. This proves Scott openness.
Saturation gives $K=\bigcap\PhiX(K)$.

Conversely, write $K=\bigcap\F$. It is saturated, and properness and the
displayed equality imply $K\ne\varnothing$. Given an open cover of $K$,
its union belongs to $\F$. The finite unions of members of the cover form
a directed family, so Scott openness gives a finite union in $\F$.
That finite union contains $K$.
\end{proof}

We recall the following minimal closed set lemma
\cite[Lemma~2.12]{SSXZ}, including its proof for completeness.

\begin{lemma}[Minimal closed set lemma]\label{lem:minimal}
Let $X$ be a $T_0$-space and $\F$ a proper Scott-open filter on $\Om(X)$, and let $C$ be a
closed subset of $X$ meeting every member of $\F$. There exists a closed
set $A\subseteq C$ minimal among the closed sets meeting every member of
$\F$. Every such minimal set is irreducible.
\end{lemma}

\begin{proof}
Consider the closed subsets of $C$ meeting every member of $\F$.
For a nonempty descending chain $\{C_i:i\in I\}$ in this collection,
suppose that $\bigcap_iC_i$ misses some $U\in\F$. Then
\[
 U\subseteq\bigcup_{i\in I}(X\setminus C_i).
\]
The union belongs to $\F$, and the open complements form a directed
family. Scott openness gives $X\setminus C_i\in\F$ for some $i$, a
contradiction. Zorn's lemma therefore gives a minimal member $A$.

If $A\subseteq C_1\cup C_2$ for closed sets $C_1,C_2$ and neither
contains $A$, both $A\cap C_1$ and $A\cap C_2$ are proper closed subsets
of $A$. There are $U_1,U_2\in\F$ with
$A\cap C_i\cap U_i=\varnothing$. But then
$A\cap U_1\cap U_2=\varnothing$, contradicting $U_1\cap U_2\in\F$.
\end{proof}

\subsection*{Irreducible sets and families}

The \emph{Smyth power space} $P_S(X)$ is the set $\K(X)$ equipped with the
upper Vietoris topology, whose basic open sets are
\[
 \Box U=\{G\in\K(X):G\subseteq U\},~ U\in\Om(X).
\]
Its specialization order is reverse inclusion.
For a nonempty family $\mathcal A\subseteq\K(X)$, extend the notation
$\PhiX$ by putting
\begin{equation}\label{eq:family-filter}
 \PhiX(\mathcal A)=\bigcup_{G\in\mathcal A}\PhiX(G)
 =\{U\in\Om(X):G\subseteq U\text{ for some }G\in\mathcal A\}.
\end{equation}
The role of irreducibility is described by the following standard fact
\cite[Lemma~2.23 and Remark~2.24]{XuH}.

\begin{lemma}\label{lem:family-filter}
Let $\varnothing\ne\mathcal A\subseteq\K(X)$. The set
$\PhiX(\mathcal A)$ is a nonempty Scott-open upper set that does not
contain $\varnothing$. It is a filter if and only if $\mathcal A$ is
irreducible in $P_S(X)$. Moreover,
\begin{equation}\label{eq:family-identities}
 \bigcap\PhiX(\mathcal A)=\bigcap\mathcal A,
 ~
 \PhiX(\mathcal A)=\PhiX\!\left(\cl{\mathcal A}^{\,P_S(X)}\right).
\end{equation}
\end{lemma}

\begin{proof}
For each $G\in\mathcal A$, Lemma~\ref{lem:compact} shows that
$\PhiX(G)$ is Scott-open in $\Om(X)$. Hence their union
$\PhiX(\mathcal A)$ is a Scott-open upper set. It contains $X$
and does not contain $\varnothing$, because $\mathcal A$ and all
its members are nonempty.

We next characterize closure under finite intersections.
Suppose first that $\mathcal A$ is irreducible in $P_S(X)$, and let
$U,V\in\PhiX(\mathcal A)$. Then $\mathcal A$ meets both
$\Box U$ and $\Box V$. By irreducibility, it meets
\[
 \Box U\cap\Box V=\Box(U\cap V).
\]
Consequently, some $G\in\mathcal A$ is contained in $U\cap V$,
so $U\cap V\in\PhiX(\mathcal A)$. Thus $\PhiX(\mathcal A)$
is a filter.

Conversely, suppose that $\PhiX(\mathcal A)$ is a filter.
Let $\mathcal V_1,\mathcal V_2$ be open subsets of $P_S(X)$ that
meet $\mathcal A$. Since the sets $\Box U$, $U\in\Om(X)$,
form a base for the upper Vietoris topology, there exist
$U_1,U_2\in\Om(X)$ such that
\[
 \mathcal A\cap\Box U_i\ne\varnothing,
 ~ \Box U_i\subseteq\mathcal V_i
 ~ (i=1,2).
\]
Thus $U_1,U_2\in\PhiX(\mathcal A)$, and the filter property gives
$U_1\cap U_2\in\PhiX(\mathcal A)$. Hence $\mathcal A$ meets
$\Box(U_1\cap U_2)\subseteq\mathcal V_1\cap\mathcal V_2$.
This proves that $\mathcal A$ is irreducible in $P_S(X)$.

Using $G=\bigcap\PhiX(G)$ for each saturated $G$, we obtain
\[
 \bigcap\PhiX(\mathcal A)
 =\bigcap_{G\in\mathcal A}\bigcap\PhiX(G)
 =\bigcap_{G\in\mathcal A}G
 =\bigcap\mathcal A.
\]
Finally, for every $U\in\Om(X)$, the open set $\Box U$ meets
$\mathcal A$ if and only if it meets
$\cl{\mathcal A}^{\,P_S(X)}$. By the definition of $\PhiX$,
this gives the second identity in \eqref{eq:family-identities}.
\end{proof}

The following combines the classical Hofmann--Mislove theorem with its
irreducible-set and Smyth-power formulations; see
\cite{HM,KP}, \cite[Theorems~3.6 and~5.7]{SSXZ}, and
\cite[Corollaries~4.3 and~4.4]{XYY}.

\begin{theorem}\label{thm:classical-forms}
For a $T_0$-space $X$, the following conditions are equivalent.
\begin{enumerate}
\item $X$ is sober.
\item Every proper Scott-open filter on $\Om(X)$ is $\PhiX(K)$ for
some $K\in\K(X)$.
\item For every irreducible subset $A\subseteq X$, there exists
$K\in\K(X)$ such that $\F_A=\PhiX(K)$.
\item For every irreducible subset $\mathcal A\subseteq P_S(X)$,
there exists $K\in\K(X)$ such that $\PhiX(\mathcal A)=\PhiX(K)$.
\end{enumerate}
In \textup{(3)} and \textup{(4)}, one may equivalently restrict the
respective test sets to irreducible closed sets. Each representation may
equivalently be written as recovery from the intersection of the filter;
the representing set is unique.
\end{theorem}

\begin{proof}
The equivalence of (1) and (2) is the classical Hofmann--Mislove theorem.
Condition (2) implies (3) and, by Lemma~\ref{lem:family-filter}, (4).
Conversely, suppose (3), fix an irreducible $A$, and write
$\F_A=\PhiX(K)$. The set $K$ meets $C=\cl A$, since otherwise the
open set $X\setminus C$ would belong to $\F_A$. Choose $x\in K\cap C$.
Every open neighbourhood of every $y\in C$ meets $A$ and hence contains
$K$, so it contains $x$. Thus $C\subseteq\cl{\{x\}}$; the reverse
inclusion follows from $x\in C$. This proves (1).

For (4)$\Rightarrow$(3), the map
$\eta_X:X\to P_S(X)$, $x\mapsto\uparrow x$, is continuous because
$\eta_X^{-1}(\Box U)=U$. Consequently $\eta_X(A)$ is irreducible and
\[
 \PhiX(\eta_X(A))=\bigcup_{a\in A}\PhiX(\uparrow a)=\F_A.
\]
Finally, $\F_A=\F_{\cl A}$ and the second identity in
\eqref{eq:family-identities} justify the closed-set restrictions.
The assertions about recovery and uniqueness follow from
Lemma~\ref{lem:compact}.
\end{proof}

\subsection*{Filtered compact families and countable open bases}

A nonempty family $\mathcal K\subseteq\K(X)$ is \emph{filtered} if
for $G,H\in\mathcal K$ some $L\in\mathcal K$ satisfies
$L\subseteq G\cap H$. Such a family is irreducible in $P_S(X)$ and
therefore induces the proper Scott-open filter in
\eqref{eq:family-filter}.
Recall that $X$ is \emph{well-filtered} if, for every such family and
every $U\in\Om(X)$,
\begin{equation}\label{eq:wf}
 \bigcap\mathcal K\subseteq U~\Longrightarrow~
 G\subseteq U\text{ for some }G\in\mathcal K.
\end{equation}
Restricting to countable filtered families gives
\emph{$\omega$-well-filteredness} \cite{XSXZ}.

\begin{proposition}[The well-filtered versions]\label{prop:wf-versions}
Let $X$ be a $T_0$-space.
\begin{enumerate}
\item $X$ is well-filtered if and only if, for every nonempty filtered
family $\mathcal K\subseteq\K(X)$, $
 \PhiX(\mathcal K)=\PhiX\!\left(\bigcap\mathcal K\right).$
 
\item $X$ is $\omega$-well-filtered if and only if the same equality
holds for every nonempty countable filtered family
$\mathcal K\subseteq\K(X)$.
\end{enumerate}
In either assertion, the equality implies that $\bigcap\mathcal K$ is
nonempty, compact and saturated. Equivalently, one may require the
existence of $K\in\K(X)$ with $\PhiX(\mathcal K)=\PhiX(K)$.
In \textup{(2)}, it suffices to test descending sequences in $\K(X)$.
\end{proposition}

\begin{proof}
For a filtered $\mathcal K$, equation~\eqref{eq:wf} says precisely that
an open set containing $\bigcap\mathcal K$ belongs to
$\PhiX(\mathcal K)$. The reverse inclusion is automatic. Apply
Lemmas~\ref{lem:compact} and~\ref{lem:family-filter} to obtain all the
representation assertions. This is the filter form of well-filteredness
in \cite[Theorem~3.11]{SSXZ}, with the same argument for its countable
version.

If $\mathcal K=\{G_n:n\in\N\}$ is filtered, choose $K_1=G_1$ and
recursively choose $K_{n+1}\in\mathcal K$ with
$K_{n+1}\subseteq K_n\cap G_{n+1}$. Then $\{K_n\}$ is a decreasing
cofinal subfamily. It has the same intersection and induces the same
filter as $\mathcal K$. This also proves the final assertion; compare
\cite[Definition~3.9 and Proposition~3.10]{XSXZ}.
\end{proof}

\begin{definition}\label{def:Hco}
For a nonempty filtered family
$\mathcal U\subseteq\Om(X)\setminus\{\varnothing\}$, put
\begin{equation}\label{eq:open-generated}
 \Psi_X(\mathcal U)
 =\bigcup_{U\in\mathcal U}\uparrow_{\Om(X)}U
 =\{V\in\Om(X):U\subseteq V\text{ for some }U\in\mathcal U\}.
\end{equation}
The space $X$ satisfies \textup{(\Hco)} if, for every such
\emph{countable} family $\mathcal U$ for which $\Psi_X(\mathcal U)$ is
Scott-open, there exists $K\in\K(X)$ such that
$\Psi_X(\mathcal U)=\PhiX(K)$.
Here ``c'' stands for countable and ``o'' for open.
\end{definition}

The Scott-open requirement in Definition~\ref{def:Hco} is an additional
hypothesis. In contrast, for filtered compact saturated families the
induced filter is automatically Scott-open. The same cofinal-sequence
argument as above allows $\mathcal U$ to be replaced by a decreasing
sequence of nonempty open sets. Hence \Hco\ is exactly representation
for proper Scott-open filters with countable open bases.

\begin{remark}[The scope of Goubault-Larrecq's theorem]\label{rem:GL-scope}
Goubault-Larrecq \cite[Theorem~1 and subsequent remark]{GL} permits
descending nonempty subsets $B_1\supseteq B_2\supseteq\cdots$, not
necessarily open, as generators:
\[
 \G=\bigcup_n\PhiX(B_n)
   =\{U\in\Om(X):B_n\subseteq U\text{ for some }n\}.
\]
Whenever $\G$ is Scott-open, his result represents it under
$\omega$-well-filteredness. Taking the $B_n$ open gives
$\omega$-well-filteredness $\Rightarrow$ \Hco.

Requiring representation for \emph{all} such Scott-open filters with
arbitrary subset generators is in fact equivalent to
$\omega$-well-filteredness: the converse follows by taking
$B_n=K_n\in\K(X)$ and applying Proposition~\ref{prop:wf-versions}.
Thus \Hco\ is specifically the open-generator case of the result.
\end{remark}

\subsection*{Irreducible subset systems and the two \texorpdfstring{$H$}{H}-versions}

We use the irreducible subset systems of \cite{XuH}; their
Hofmann--Mislove formulation is developed in \cite{SSXZ}.
Here $\mathbf{Top}_0$ is the category of $T_0$-spaces and continuous
maps, and $\mathbf{Set}$ is the category of sets and maps.

\begin{definition}\label{def:H-system}
An \emph{irreducible subset system} is a functor
$H:\mathbf{Top}_0\to\mathbf{Set}$ such that, for each space $X$,
\[
 \{\{x\}:x\in X\}\subseteq H(X)
 \subseteq\{A\subseteq X:A\text{ is nonempty and irreducible}\},
\]
and, for every continuous map $f:X\to Y$, the map $H(f)$ is given by
$H(f)(A)=f(A)$. A member of $H(X)$ is called an \emph{$H$-set}.
Put
\[
 H_c(X)=\{\cl A:A\in H(X)\}.
\]
A \emph{closed $H$-set} means a member of $H_c(X)$.
The space $X$ is \emph{$H$-sober} if every $A\in H(X)$ has
$\cl A=\cl{\{x\}}$ for a unique $x\in X$.
It is \emph{super $H$-sober} if $P_S(X)$ is $H$-sober.
\end{definition}

The use of $H_c(X)$ matters: it consists of closures of $H$-sets,
rather than just those members of $H(X)$ that happen to be closed.
The two filter classes below consist of proper Scott-open filters,
by \eqref{eq:FA} and Lemma~\ref{lem:family-filter}.
We first state the point-induced formulation of
\cite[Theorems~3.6 and~5.7]{SSXZ}.

\begin{theorem}[The point-induced $H$-version]\label{thm:H-point}
Let $H$ be an irreducible subset system. For a $T_0$-space $X$, the
following conditions are equivalent.
\begin{enumerate}
\item $X$ is $H$-sober.
\item For every $A\in H(X)$, $\F_A=\PhiX\!\left(\bigcap\F_A\right)$.
 \item For every $A\in H(X)$, the set $\bigcap\F_A$ belongs to $\K(X)$
and $\F_A=\PhiX(\bigcap\F_A)$.
\item For every $A\in H(X)$, there exists $K\in\K(X)$ such that
$\F_A=\PhiX(K)$.
\end{enumerate}
The representing set is unique and equals
$\bigcap\F_A=\uparrow x$, where $\cl A=\cl{\{x\}}$.
The same equivalences hold with the test sets taken from $H_c(X)$.
\end{theorem}

\begin{proof}
If $X$ is $H$-sober, equation~\eqref{eq:point-filter} gives (4).
Conversely, suppose (4), fix $A\in H(X)$, and write
$\F_A=\PhiX(K)$. Since $X\setminus\cl A\notin\F_A$, there is
$x\in K\cap\cl A$. Every open neighbourhood of any $y\in\cl A$ meets
$A$, belongs to $\F_A$, and therefore contains $x$. Thus
$\cl A\subseteq\cl{\{x\}}$. The reverse inclusion follows from
$x\in\cl A$, and uniqueness follows from $T_0$.
This proves (1).

Lemma~\ref{lem:compact} gives the equivalence of (2), (3) and (4),
as well as uniqueness of the representing set.
Equation~\eqref{eq:point-filter} identifies it as $\uparrow x$.
Finally, $\F_A=\F_{\cl A}$ gives the formulation using $H_c(X)$.
\end{proof}

The second formulation concerns $H$-sets of the Smyth power space.
It is the filter representation form of
\cite[Theorem~5.5]{XuH}.

\begin{theorem}[The Smyth $H$-version]\label{thm:H-Smyth}
Let $H$ be an irreducible subset system. For a $T_0$-space $X$, the
following conditions are equivalent.
\begin{enumerate}
\item $P_S(X)$ is $H$-sober, that is, $X$ is super $H$-sober.
\item For every $\mathcal A\in H(P_S(X))$, $\PhiX(\mathcal A)=\PhiX\!\left(\bigcap\mathcal A\right)$.
\item For every $\mathcal A\in H(P_S(X))$, the set
$\bigcap\mathcal A$ belongs to $\K(X)$ and
$\PhiX(\mathcal A)=\PhiX(\bigcap\mathcal A)$.
\item For every $\mathcal A\in H(P_S(X))$, there exists $K\in\K(X)$
such that $\PhiX(\mathcal A)=\PhiX(K)$.
\end{enumerate}
The representing set is unique and equals $\bigcap\mathcal A$.
Equivalently, the test families may be taken from $H_c(P_S(X))$.
\end{theorem}

\begin{proof}
For $\mathcal A\in H(P_S(X))$, $K\in\K(X)$ and $U\in\Om(X)$,
\[
 \mathcal A\cap\Box U\ne\varnothing
 \ \Longleftrightarrow\ U\in\PhiX(\mathcal A);
 ~~
 K\in\Box U\ \Longleftrightarrow\ U\in\PhiX(K).
\]
Since the sets $\Box U$ form a base, it follows that
\begin{equation}\label{eq:Smyth-generic}
 \PhiX(\mathcal A)=\PhiX(K)
 ~\Longleftrightarrow~
 \cl{\mathcal A}^{\,P_S(X)}=\cl{\{K\}}^{\,P_S(X)}.
\end{equation}
This proves (1)$\Leftrightarrow$(4).
The filter $\PhiX(\mathcal A)$ is proper and Scott-open, and
$\bigcap\PhiX(\mathcal A)=\bigcap\mathcal A$, by
Lemma~\ref{lem:family-filter}. Lemma~\ref{lem:compact} therefore gives
the remaining equivalences and uniqueness.
The second identity in \eqref{eq:family-identities} gives the
formulation using $H_c(P_S(X))$.
\end{proof}

These two classes are related by the continuous map
\[
 \eta_X:X\longrightarrow P_S(X),~ x\longmapsto\uparrow x,
 ~ \eta_X^{-1}(\Box U)=U.
\]
Functoriality gives $\eta_X(A)\in H(P_S(X))$ for $A\in H(X)$, and
$\PhiX(\eta_X(A))=\F_A$. Consequently,
\begin{equation}\label{eq:H-filter-inclusion}
 \{\F_A:A\in H(X)\}
 \subseteq
 \{\PhiX(\mathcal A):\mathcal A\in H(P_S(X))\}.
\end{equation}
In particular, super $H$-sobriety always implies $H$-sobriety.
The class on the right also contains every $\PhiX(K)$ with
$K\in\K(X)$, by taking $\mathcal A=\{K\}$.

\subsection*{Property \texorpdfstring{$Q$}{Q} and countable restrictions}

Property $Q$ was introduced in \cite[Definition~5.10]{XuH}.
We use its local form from \cite[Definition~6.2]{SSXZ}.

\begin{definition}\label{def:Q}
An irreducible subset system $H$ has \emph{property $Q$ for $X$} if,
whenever $\mathcal A\in H(P_S(X))$ and a closed set $C\subseteq X$
meets every $G\in\mathcal A$, there exists $B\in H_c(X)$ such that
\[
 B\subseteq C,~ B\cap G\ne\varnothing
 ~\text{for every }G\in\mathcal A.
\]
It has \emph{property $Q$} if it has property $Q$ for every $T_0$-space.
\end{definition}

The following exact comparison is
\cite[Theorem~6.4]{SSXZ}. We include a proof to exhibit how the two
filter representations enter the argument.

\begin{theorem}[The property-$Q$ comparison]\label{thm:H-Q}
Let $H$ be an irreducible subset system and $X$ a $T_0$-space. Then
\[
 P_S(X)\text{ is }H\text{-sober}
 ~\Longleftrightarrow~
 \left\{
 \begin{array}{l}
 X\text{ is }H\text{-sober},\\
 H\text{ has property }Q\text{ for }X.
 \end{array}\right.
\]
In particular, if $H$ has property $Q$, the representation properties
in Theorems~\ref{thm:H-point} and~\ref{thm:H-Smyth} are equivalent
for every $T_0$-space.
\end{theorem}

\begin{proof}
Suppose that $P_S(X)$ is $H$-sober. Equation~\eqref{eq:H-filter-inclusion}
and the two representation theorems show that $X$ is $H$-sober.
Let $\mathcal A$ and $C$ be as in Definition~\ref{def:Q}.
Theorem~\ref{thm:H-Smyth} gives
$K=\bigcap\mathcal A\in\K(X)$ and
$\PhiX(\mathcal A)=\PhiX(K)$. If $C\cap K=\varnothing$, then
$X\setminus C\in\PhiX(\mathcal A)$, so some $G\in\mathcal A$
misses $C$, a contradiction.
Choose $x\in C\cap K$. The closed set
$B=\cl{\{x\}}\in H_c(X)$ is contained in $C$ and meets every
$G\in\mathcal A$, proving property $Q$ for $X$.

Conversely, suppose that $X$ is $H$-sober and $H$ has property $Q$
for $X$. Let $\mathcal A\in H(P_S(X))$ and suppose that
$\bigcap\mathcal A\subseteq U$ for an open set $U$.
If no $G\in\mathcal A$ is contained in $U$, the closed set
$C=X\setminus U$ meets every $G\in\mathcal A$.
Property $Q$ gives $B\in H_c(X)$, contained in $C$, that still meets
every $G$. By $H$-sobriety, $B=\cl{\{x\}}$ for some $x\in X$.
Each saturated $G$ meeting $\cl{\{x\}}$ contains $x$.
Hence $x\in C\cap\bigcap\mathcal A$, a contradiction.
We have proved
$\PhiX(\bigcap\mathcal A)\subseteq\PhiX(\mathcal A)$; the opposite
inclusion is automatic. Apply Theorem~\ref{thm:H-Smyth}.
\end{proof}

\begin{definition}[Countable restriction]\label{def:H-countable}
For an irreducible subset system $H$, define
\[
 H^\omega(X)=\{A\in H(X):A\text{ is countable}\}.
\]
This is again an irreducible subset system, since direct images preserve
countability. Its closed-set family is
\[
 (H^\omega)_c(X)
 =\{\cl A:A\in H(X),\ A\text{ is countable}\}.
\]
In particular, writing $\mathcal R(X)$ for all nonempty irreducible
subsets of $X$, its countable restriction is
\[
 \mathcal R^\omega(X)
 =\{A\subseteq X:A\text{ is countable, nonempty and irreducible}\}.
\]
\end{definition}

Thus $(H^\omega)_c(X)$ consists of closures of countable $H$-sets;
its members need not themselves be countable. In particular,
\begin{equation}\label{eq:omega-Romega}
 X\text{ is }\omega\text{-sober}
 ~\Longleftrightarrow~
 X\text{ is }\mathcal R^\omega\text{-sober}.
\end{equation}
Applying Theorems~\ref{thm:H-point} and~\ref{thm:H-Smyth} to
$H^\omega$ gives the two countable $H$-versions: the point-induced
filters characterize $H^\omega$-sobriety of $X$, and the
Smyth-induced filters characterize $H^\omega$-sobriety of $P_S(X)$.
Theorem~\ref{thm:H-Q} compares them when $H^\omega$ has property $Q$
for $X$. For $\mathcal R^\omega$ this property will be proved in
Lemma~\ref{lem:Romega-Q}.

For a countable family $\mathcal A\subseteq P_S(X)$, the
countability condition concerns the number of compact saturated members
of $\mathcal A$. Each such member may be uncountable, and the induced
Scott-open filter need not have a countable open base.

\section{An enlarged class of Scott-open filters}

\begin{definition}\label{def:cb}
A proper filter $\F\subseteq\Om(X)$ is \emph{countably generated},
or \emph{countably based}, if there are $U_n\in\F$, $n\in\N$, such that
\begin{equation}\label{eq:base}
 \F=\{U\in\Om(X):U_n\subseteq U\text{ for some }n\in\N\}.
\end{equation}
The base may equivalently be required to be a countable filtered family
of nonempty open sets. Indeed, finite intersections give a decreasing
base; conversely, from a countable filtered open base one can choose a
decreasing cofinal sequence of its members. Thus we may assume
$U_{n+1}\subseteq U_n$. These sets form a base of the filter, not
necessarily a base of the topology of $X$.

Let $\OFcb(X)$ denote the family of proper Scott-open filters on
$\Om(X)$ having a countable open base, where the subscript ``cb''
stands for ``countably based''. We define the enlarged class
\begin{equation}\label{eq:enlarged}
 \Cw(X)=\OFcb(X)\cup
 \{\F_A:A\subseteq X\text{ is countable and irreducible}\}.
\end{equation}
Every member of $\Cw(X)$ is a proper Scott-open filter. Members of the
second subfamily need not have countable open bases.
\end{definition}

By Theorem~\ref{thm:H-point}, applied to $H=\mathcal R^\omega$,
$\omega$-sobriety is already equivalent to compact saturated
representation of all filters $\F_A$ induced by countable
irreducible subsets $A$ of $X$. The next theorem shows that
adjoining $\OFcb(X)$ to this test class preserves the equivalence.
Its additional assertion is that $\omega$-sobriety also represents
every proper Scott-open filter having a countable open base.
In general, the filters $\F_A$ cannot be omitted: representation
for $\OFcb(X)$ alone is strictly weaker, as
Remark~\ref{rem:scope} and Example~\ref{ex:cofinite} explain.

\begin{theorem}[A countable version of the Hofmann--Mislove theorem]\label{thm:cb}
For a $T_0$-space $X$, the following conditions are equivalent:
\begin{enumerate}[label=\textup{(\arabic*)}]
\item $X$ is $\omega$-sober.
\item For every $\F\in\Cw(X)$, $
 \F=\PhiX\!\left(\bigcap\F\right)$.
\item For every $\F\in\Cw(X)$, the set $\bigcap\F$ is nonempty,
compact and saturated, and
$
 \F=\PhiX\!\left(\bigcap\F\right)$.
\item For every $\F\in\Cw(X)$, there exists a nonempty compact
saturated subset $K$ of $X$ such that $\F=\PhiX(K)$.
\end{enumerate}
The representing set in \textup{(4)} is unique and equals $\bigcap\F$.
For $\F=\F_A$, it is $\uparrow x$, where $x$ is the unique point
satisfying $\cl A=\cl{\{x\}}$.
\end{theorem}

\begin{proof}
We first prove (1)$\Rightarrow$(4). Suppose that $X$ is $\omega$-sober
and $\F\in\OFcb(X)$, with a countable open base $\{U_n:n\in\N\}$.
Let $W\in\Om(X)\setminus\F$. Every member of $\F$
meets $C=X\setminus W$: otherwise it would be contained in $W$, forcing
$W\in\F$. Apply Lemma~\ref{lem:minimal} to obtain a minimal irreducible
closed set $A\subseteq C$ meeting every member of $\F$.

Choose $x_n\in A\cap U_n$ and put $D=\{x_n:n\in\N\}$. The closed
set $\cl D\subseteq A$ still meets every member of $\F$, because every
member contains some $U_n$. Minimality gives $A=\cl D$.
Thus $A$ is an irreducible closed set with a countable dense subset.
By $\omega$-sobriety, $A=\cl{\{x\}}$ for some $x\in X$. Every
$U\in\F$ meets $A$, so it contains $x$. Since $x\in A\subseteq C$,
$
 x\in\bigcap\F\setminus W$.
Consequently an open set contains $\bigcap\F$ if and only if it belongs
to $\F$. Lemma~\ref{lem:compact} supplies the required nonempty compact
saturated representing set.

For the other subfamily in \eqref{eq:enlarged}, let $A$ be countable
and irreducible. By $\omega$-sobriety, $\cl A=\cl{\{x\}}$ for some
$x\in X$. Equation~\eqref{eq:point-filter} gives
$\F_A=\PhiX(\uparrow x)$, proving (4) for every member of $\Cw(X)$.

The implication (4)$\Rightarrow$(2) follows from
$K=\bigcap\PhiX(K)$ for saturated $K$. The implication
(2)$\Rightarrow$(3) is Lemma~\ref{lem:compact}, and
(3)$\Rightarrow$(4) is immediate.

Finally, (2) applied to the filters $\F_A$ induced by countable
irreducible subsets is condition~(2) of Theorem~\ref{thm:H-point}
for $H=\mathcal R^\omega$. That theorem gives $\omega$-sobriety
and hence (1). Uniqueness of the representing set follows from
saturation, and \eqref{eq:point-filter} gives its stated form for $\F_A$.
\end{proof}

A subset $K$ has \emph{countable outer character} if its open-neighbourhood
filter has a countable open base. Let $\K_\cb(X)$ consist of the members
of $\K(X)$ having this property.

\begin{corollary}\label{cor:iso}
If $X$ is $\omega$-sober, then
\[
 \PhiX:(\K_\cb(X),\supseteq)
 \longrightarrow(\OFcb(X),\subseteq)
\]
is an order isomorphism, with inverse $\F\mapsto\bigcap\F$.
With inclusion on both sides it is an order anti-isomorphism.
\end{corollary}

\begin{proof}
Theorem~\ref{thm:cb} gives the bijection and its inverse.
For saturated $K,L$, the equivalence
$K\supseteq L\Longleftrightarrow\PhiX(K)\subseteq\PhiX(L)$ follows
by intersecting their open neighbourhoods.
\end{proof}

\begin{remark}\label{rem:scope}
The restriction to $\K_\cb(X)$ is necessary: a compact saturated set
need not have a countable open-neighbourhood base.

The enlargement in \eqref{eq:enlarged} is necessary in general for
the converse in Theorem~\ref{thm:cb}. On an uncountable set $X$
with the cofinite topology, every proper countably generated filter
has the form
\[
 \PhiX(X\setminus E),~ E\subseteq X\text{ countable},
\]
and $X\setminus E$ is nonempty, compact and saturated.
Consequently, all of conditions~\textup{(2)--(4)} in
Theorem~\ref{thm:cb} hold when their test class is restricted to
$\OFcb(X)$. Nevertheless, $X$ is not $\omega$-sober. For every
countably infinite $D\subseteq X$, the filter
$
 \F_D=\Om(X)\setminus\{\varnothing\}$ 
has empty intersection and no countable open base. It belongs to
$\Cw(X)\setminus\OFcb(X)$ and fails the representation condition.
The details are given in Example~\ref{ex:cofinite}.

If $X$ is first-countable, however, then
$\Cw(X)=\OFcb(X)$. Indeed, for a countable irreducible subset
$A\subseteq X$, choose a countable base $\mathcal B_a$ of open
neighbourhoods at each $a\in A$. The countable family
$\bigcup_{a\in A}\mathcal B_a$ is an open base of $\F_A$.
Thus, within the class of first-countable $T_0$-spaces,
representation for $\OFcb(X)$ alone does characterize
$\omega$-sobriety.

We retain the usual meaning of \emph{countably generated}; the
whole class $\Cw(X)$ need not consist of countably generated filters.
\end{remark}

\subsection*{The countable irreducible subset system}

We now connect the enlarged-filter theorem with both general
$H$-versions. The following countable selection argument is the
essential step.

\begin{lemma}\label{lem:Romega-Q}
The countable irreducible subset system $\mathcal R^\omega$ has
property $Q$.
\end{lemma}

\begin{proof}
Let $X$ be a $T_0$-space,
$\mathcal A=\{K_n:n\in\N\}\in\mathcal R^\omega(P_S(X))$, and
let $C$ be a closed set meeting every $K_n$.
For a finite family, we allow repetitions in this enumeration.
By Lemma~\ref{lem:family-filter},
\[
 \F=\PhiX(\mathcal A)=\bigcup_{n\in\N}\PhiX(K_n)
\]
is a proper Scott-open filter.
A closed set $B$ meets every member of $\F$ if and only if it meets
every $K_n$. Indeed, meeting each $K_n$ suffices, since each member
of $\F$ contains some $K_n$. Conversely, if $B$ misses $K_n$, then
$X\setminus B\in\PhiX(K_n)\subseteq\F$ witnesses failure to meet
every member of $\F$.

Apply Lemma~\ref{lem:minimal} to obtain a minimal irreducible closed
set $B\subseteq C$ meeting every $K_n$. Choose $x_n\in B\cap K_n$
and put $D=\{x_n:n\in\N\}$. The closed set $\cl D\subseteq B$
still meets every $K_n$, so minimality gives
$
 B=\cl D$.
Since $B$ is irreducible, $D$ is irreducible. Thus
$D\in\mathcal R^\omega(X)$ and
$B\in(\mathcal R^\omega)_c(X)$, as required by
Definition~\ref{def:Q}.
\end{proof}

The lemma requires the family $\mathcal A$ to be countable and
irreducible in $P_S(X)$; it need not be filtered under inclusion.
It proves that the closed witness has a countable dense irreducible
subset, without requiring the closed witness itself to be countable.

\begin{theorem}[The countable Smyth version]\label{thm:omega-Smyth}
For a $T_0$-space $X$, the following conditions are equivalent.
\begin{enumerate}
\item $X$ is $\omega$-sober, that is, $\mathcal R^\omega$-sober.
\item For every countable irreducible subset $A\subseteq X$, there
exists $K\in\K(X)$ such that $\F_A=\PhiX(K)$.
\item $P_S(X)$ is $\omega$-sober.
\item For every countable irreducible subset
$\mathcal A\subseteq P_S(X)$, $\PhiX(\mathcal A)=\PhiX\!\left(\bigcap\mathcal A\right)$.
\item For every $\mathcal A$ as in \textup{(4)},
$\bigcap\mathcal A$ is nonempty, compact and saturated, and
$\PhiX(\mathcal A)=\PhiX(\bigcap\mathcal A)$.
\item For every $\mathcal A$ as in \textup{(4)}, there exists
$K\in\K(X)$ such that $\PhiX(\mathcal A)=\PhiX(K)$.
\end{enumerate}
These conditions are also equivalent to the conditions of
Theorem~\ref{thm:cb}. All representing compact saturated sets are
unique: in \textup{(2)} the set is $\uparrow x$, where
$\cl A=\cl{\{x\}}$, and in \textup{(6)} it is
$\bigcap\mathcal A$.

In \textup{(2)}, one may instead test irreducible closed subsets
of $X$ having a countable dense subset. In \textup{(4)--(6)},
one may instead test irreducible closed subsets of $P_S(X)$
having a countable dense subfamily.
\end{theorem}

\begin{proof}
Equation~\eqref{eq:omega-Romega} and Theorem~\ref{thm:H-point}
give (1)$\Leftrightarrow$(2).
Lemma~\ref{lem:Romega-Q} and Theorem~\ref{thm:H-Q}, applied to
$H=\mathcal R^\omega$, give (1)$\Leftrightarrow$(3).
Theorem~\ref{thm:H-Smyth} gives the equivalence of (3)--(6)
and identifies the representing set in (6).
Theorem~\ref{thm:H-point} identifies the set in (2).
The closed-set formulations follow from the descriptions of
$(\mathcal R^\omega)_c(X)$ and
$(\mathcal R^\omega)_c(P_S(X))$.
\end{proof}

\section{Traces on countable subspaces}

A further equivalent formulation uses filters pulled back from countable
subspaces. For $D\subseteq X$ let
\[
 \rho_D:\Om(X)\longrightarrow\Om(D),~ \rho_D(U)=U\cap D,
\]
where $D$ carries the subspace topology. This map preserves finite
intersections and arbitrary unions. Hence, if $\G$ is a proper Scott-open
filter on $\Om(D)$, then
\begin{equation}\label{eq:pullback}
 \rho_D^{-1}(\G)=\{U\in\Om(X):U\cap D\in\G\}
\end{equation}
is a proper Scott-open filter on $\Om(X)$. Its membership is determined
by the traces of open sets on $D$.

For irreducible $D$, the completely prime filter $\F_D$ from
\eqref{eq:FA} is the pullback of $\Om(D)\setminus\{\varnothing\}$.

\begin{theorem}\label{thm:exact}
For a $T_0$-space $X$, the following conditions are equivalent.
\begin{enumerate}
\item $X$ is $\omega$-sober.
\item For every countable subspace $D\subseteq X$ and every proper
Scott-open filter $\G\subseteq\Om(D)$, there exists $K\in\K(X)$ with $
 \rho_D^{-1}(\G)=\PhiX(K)$.

\item For every $D$ and $\G$ as in (2), writing
$\F=\rho_D^{-1}(\G)$, one has
$
 \F=\PhiX\!\left(\bigcap\F\right)$.

\item For every countable irreducible subset $D\subseteq X$,
$
 \F_D=\PhiX\!\left(\bigcap\F_D\right)$.
\item For every countable irreducible subset $D\subseteq X$,
$
 \cl D\cap\bigcap\F_D\ne\varnothing$.
\end{enumerate}
In (2) the representing compact saturated set is unique. In (4), it
is $\uparrow x$, where $x$ is the unique point with
$\cl D=\cl{\{x\}}$.
\end{theorem}

\begin{proof}
Assume (1), and fix a countable $D$, a proper Scott-open $\G$, and
$\F=\rho_D^{-1}(\G)$. Let $C$ be a closed set meeting every member of
$\F$. By Lemma~\ref{lem:minimal}, it contains a minimal irreducible
closed set $A$ with that property.

We claim that $A\cap D$ still meets every member of $\F$. Otherwise,
for some $U\in\F$ we would have
\[
 U\cap D\subseteq (X\setminus A)\cap D.
\]
The left side belongs to $\G$, so upperness of $\G$ would give
$X\setminus A\in\F$, a contradiction. Thus $\cl{A\cap D}\subseteq A$
is a closed set meeting every member of $\F$. By minimality,
$
 A=\cl{A\cap D}$.
The set $A\cap D$ is countable. Consequently $A=\cl{\{x\}}$ for some
$x\in X$, and $x\in C\cap\bigcap\F$.

We have proved the closed-set intersection property for $\F$ and every
closed $C$. Its contrapositive, applied to complements of open sets,
gives $\F=\PhiX(\bigcap\F)$. Lemma~\ref{lem:compact} proves (2).

The implication (2)$\Rightarrow$(3) follows by saturation. For
(3)$\Rightarrow$(4), take a countable irreducible subspace $D$ and
$\G=\Om(D)\setminus\{\varnothing\}$. This is a proper completely prime
filter, and its inverse image is $\F_D$.

Assume (4). The open set $X\setminus\cl D$ does not belong to $\F_D$,
so it cannot contain $\bigcap\F_D$. This proves (5).

Finally, assume (5) and choose $x\in\cl D\cap\bigcap\F_D$. Every
open set meeting $D$ contains $x$. If $y\in\cl D$, each open neighbourhood
of $y$ meets $D$ and therefore contains $x$. Hence
$y\in\cl{\{x\}}$, proving $\cl D\subseteq\cl{\{x\}}$. The reverse
inclusion follows from $x\in\cl D$. Uniqueness follows from $T_0$.
Moreover, $\F_D=\{U:x\in U\}=\PhiX(\uparrow x)$.
\end{proof}

\begin{remark}\label{rem:two-countabilities}
The filters in Theorem~\ref{thm:exact} are not required to have countable
open bases. For example,
\[
 \F_D=\bigcup_{d\in D}\{U\in\Om(X):d\in U\}
\]
is a countable union of point-neighbourhood filters when $D$ is countable.
This does not make it a countably generated filter on $\Om(X)$.
The filter of all nonempty opens in Example~\ref{ex:cofinite} is an
explicit instance. The completely prime formulation (4) is the
$H=\mathcal R^\omega$ case of the general approach in
\cite[Theorems~3.6 and~5.7]{SSXZ}.
\end{remark}

\section{\texorpdfstring{$\omega$}{omega}-well-filteredness and strictness}

We now compare $\omega$-sobriety with the filtered-family and
open-generator representation properties from Section~2.
Write $\mathcal D(X)$ for the nonempty directed subsets of $X$ in
the specialization order, and let $\mathcal D^\omega$ be its countable
restriction. Both are irreducible subset systems.
The $\mathcal D$-sobriety condition is the usual $d$-space property.
Its countable version, the $\omega$-$d$-space property, says that the
closure of every countable directed set is a point closure
\cite{XuH,XSXZ}.

Since the specialization order of $P_S(X)$ is reverse inclusion,
its directed subsets are exactly the filtered families in $\K(X)$.
Thus Proposition~\ref{prop:wf-versions} is the specialization of
Theorem~\ref{thm:H-Smyth} to $\mathcal D$ and
$\mathcal D^\omega$, respectively. Together with
Theorems~\ref{thm:classical-forms} and~\ref{thm:omega-Smyth},
this gives Table~\ref{tab:H-versions}.

\begin{table}[htb]
\centering
\caption{Specializations of the two $H$-representation properties.}
\label{tab:H-versions}
\small
\begin{tabular*}{\textwidth}{@{\extracolsep{\fill}}lcc@{}}
\toprule
$H$ & Point-induced version & Smyth version \\
& ($X$ is $H$-sober) & ($P_S(X)$ is $H$-sober) \\
\midrule
$\mathcal R$ & Sobriety & Sobriety \\
$\mathcal D$ & $d$-space property & Well-filteredness \\
$\mathcal D^\omega$ & $\omega$-$d$-space property
  & $\omega$-well-filteredness \\
$\mathcal R^\omega$ & $\omega$-sobriety & $\omega$-sobriety \\
\bottomrule
\end{tabular*}
\end{table}

The distinction between $\omega$-well-filteredness and
$\omega$-sobriety is therefore the distinction between testing
countable \emph{filtered} compact families and testing all countable
\emph{irreducible} subsets of $P_S(X)$.

\begin{proposition}\label{prop:omega-wf}
Every $\omega$-sober space is $\omega$-well-filtered.
\end{proposition}

\begin{proof}
Every nonempty countable filtered family in $\K(X)$ is a countable
irreducible subset of $P_S(X)$. If $X$ is $\omega$-sober,
Theorem~\ref{thm:omega-Smyth}(4) gives the recovery equality for
each such family. Proposition~\ref{prop:wf-versions}(2) now yields
$\omega$-well-filteredness.
\end{proof}

With \Hco\ as in Definition~\ref{def:Hco}, Theorem~\ref{thm:cb} and
Proposition~\ref{prop:omega-wf}, together with
the known representation theorem under $\omega$-well-filteredness
\cite[Theorem~1 and subsequent remark]{GL}, give
\begin{equation}\label{eq:chain}
 \omega\text{-sobriety}
 ~\Longrightarrow~\omega\text{-well-filteredness}
 ~\Longrightarrow~\Hco.
\end{equation}
By Theorem~\ref{thm:cb}, the first property is equivalent to the
representation property for all filters in $\Cw(X)$. Both implications
in \eqref{eq:chain} are strict, as the following two examples show.

\begin{example}[Well-filtered but not $\omega$-sober]\label{ex:ultrafilter}
Let $\mathscr U$ be a free ultrafilter on $\N$ and equip $X=\N$ with
\[
 \tau=\{\varnothing\}\cup\mathscr U.
\]
This is a topology. Every cofinite set belongs to $\mathscr U$, so $X$ is
$T_1$. Two nonempty open sets have nonempty intersection, so $X$ is an
irreducible closed set. It is countable and has no generic point, because
point closures are singletons. Thus $X$ is not $\omega$-sober.

We prove that its compact subsets are exactly the finite subsets.
Suppose $K\subseteq\N$ is infinite. Partition $K$ into two infinite
pieces; at least one piece, say $B$, does not belong to $\mathscr U$.
Enumerate $B=\{b_1,b_2,\ldots\}$ and put
\[
 V_n=(\N\setminus B)\cup\{b_1,\ldots,b_n\}.
\]
Each $V_n$ is open, since $\N\setminus B\in\mathscr U$. These sets
increase to $X$, but no finite subfamily covers $K$. Hence $K$ is not
compact.

Since $X$ is $T_1$, its nonempty compact saturated sets are exactly the
nonempty finite sets. In any filtered family of such sets, choose a member
$K_0$ of minimal cardinality. A common refinement of $K_0$ and any other
member must equal $K_0$. Thus $K_0$ is contained in every member, and the
intersection of the family is $K_0$. This proves well-filteredness,
and hence $\omega$-well-filteredness. The example uses the usual existence
of a free ultrafilter in ZFC.
\end{example}

\begin{example}[The uncountable cofinite space]\label{ex:cofinite}
Let $X$ be an uncountable set with the cofinite topology. Every subset of
$X$ is compact, since its subspace topology is cofinite. Every subset is
also saturated, since $X$ is $T_1$.

For a proper countably generated filter $\F$, choose a decreasing open
base $U_n=X\setminus E_n$, where each $E_n$ is finite and
$E_n\subseteq E_{n+1}$. Put
\[
 E=\bigcup_nE_n,~ K=X\setminus E.
\]
The set $E$ is countable; hence $K$ is nonempty, compact and saturated.
It is plainly $\bigcap\F$. If an open set $W$ contains $K$, its finite
complement is contained in $E$, and therefore in some $E_n$. Thus
$U_n\subseteq W$ and $W\in\F$. We obtain
$
 \F=\PhiX(K)$.
In fact this argument represents every proper countably generated
filter, without initially assuming Scott openness. Conversely, every
countable $E\subseteq X$ is the union of an increasing sequence of
finite subsets $E_n$, and $\{X\setminus E_n:n\in\N\}$ is an open
base of $\PhiX(X\setminus E)$. Lemma~\ref{lem:compact} therefore gives
\[
 \OFcb(X)
 =\{\PhiX(X\setminus E):E\subseteq X\text{ is countable}\}.
\]
In particular, $X$ satisfies \textup{(\Hco)}, and all of
conditions~\textup{(2)--(4)} in Theorem~\ref{thm:cb} hold if the
test class $\Cw(X)$ is replaced by $\OFcb(X)$.

Choose distinct $d_n\in X$ and let
$K_n=\{d_m:m\geq n\}$. These are nonempty compact saturated sets,
they decrease, and $\bigcap_nK_n=\varnothing$. Equation~\eqref{eq:wf}
fails for $U=\varnothing$, so $X$ is not $\omega$-well-filtered.

The countable set $D=\{d_n:n\in\N\}$ is irreducible and dense in $X$.
Its closure $X$ has no generic point, so $X$ is not $\omega$-sober either.
Here
\[
 \F_D=\Om(X)\setminus\{\varnothing\}
\]
is a completely prime Scott-open filter whose intersection is empty.
In fact,
\[
 \PhiX(\{K_n:n\in\N\})=\F_D,
\]
because every nonempty cofinite open set contains some tail $K_n$.
Thus this filter is induced by a countable decreasing family of compact
saturated sets, although it has no countable open base.
Indeed, for any proposed countable open base
$X\setminus E_n$, a point $p\notin\bigcup_nE_n$ belongs to every base
member, whereas the nonempty open set $X\setminus\{p\}$ contains none
of them. Thus the countability of $D$ does not imply the countable
generation of $\F_D$. In particular, $\F_D\in\Cw(X)$ witnesses
failure of the enlarged representation property:
\[
 \F_D\ne\PhiX\!\left(\bigcap\F_D\right)
       =\PhiX(\varnothing)=\Om(X).
\]
Thus the equivalent conditions in Theorem~\ref{thm:cb} are strictly
stronger than \textup{(\Hco)}.

This example also separates the two general $H$-versions. Since $X$
is $T_1$, every directed subset is a singleton, so $X$ is both
$\mathcal D$-sober and $\mathcal D^\omega$-sober. Its failure of
$\omega$-well-filteredness shows that it is neither super
$\mathcal D$-sober nor super $\mathcal D^\omega$-sober.
By Theorem~\ref{thm:H-Q}, neither $\mathcal D$ nor
$\mathcal D^\omega$ has property $Q$ for this $X$.
\end{example}

\begin{remark}[Strict weakening of sobriety]
For completeness, $\omega$-sobriety is itself weaker than sobriety.
On an uncountable set with the cocountable topology, each countable
subspace is closed and discrete. Consequently, every countable irreducible
subset is a singleton, so the space is $\omega$-sober. However, the whole space is
irreducible and has no generic point. It is therefore not sober.
\end{remark}

\section*{Scope of the countable version}

The two natural classes of filters induced by an irreducible subset
system $H$ encode $H$-sobriety of $X$ and of $P_S(X)$.
Their relation is given by the local property-$Q$ criterion
\cite[Theorem~6.4]{SSXZ}. The point-induced formulation belongs to
the framework of \cite[Theorems~3.6 and~5.7]{SSXZ}, and the Smyth
formulation is the filter version of \cite[Theorem~5.5]{XuH}.
Theorem~\ref{thm:omega-Smyth} is the countable specialization obtained
by proving that $\mathcal R^\omega$ has property $Q$.
In particular, the $\omega$-sober spaces considered here are exactly
the $\mathcal R^\omega$-sober spaces.

Theorem~\ref{thm:cb} incorporates the proper Scott-open filters with
countable open bases into an exact representation characterization,
by adjoining the filters induced by countable irreducible subsets.
Theorem~\ref{thm:exact} supplies the equivalent countable-trace
formulation. The Smyth formulation tests countable irreducible
families of compact saturated sets; restricting those families to
filtered ones gives $\omega$-well-filteredness.

The examples distinguish these conditions from the open-generator
property \Hco. The implication
$\omega$-well-filteredness $\Rightarrow$ \Hco\ is the open-generator
case of \cite{GL}, and both implications in \eqref{eq:chain} are
strict. These distinctions keep the different roles of countability
explicit throughout the representation theorems.


\section*{Declaration of generative AI and AI-assisted technologies
in the manuscript preparation process}
During the revision of this manuscript, ChatGPT (OpenAI) was used
to assist with checking mathematical arguments, clarifying definitions
and proofs, organizing the exposition, and language editing.
The authors take responsibility for the content of the manuscript.

\end{document}